\documentclass{article}
\usepackage{amsmath,amssymb}
\newtheorem{df}{Definition}[section]
\newtheorem{thm}[df]{Theorem}
\newtheorem{lem}[df]{Lemma}
\title{On the uncertainty product of spherical wavelets}
\author{Ilona Iglewska-Nowak\footnote{West Pomeranian University of Technology in Szczecin, School of Mathematics, al. Pias\-t\'ow 17, 70--310 Szczecin, Poland}}

\begin{document}

\maketitle

\bibliographystyle{amsplain}

\begin{abstract}In the paper, asymptotic behavior of the uncertainty product for a family of zonal spherical wavelets is computed. The family contains the most popular wavelets, such as Gauss--Weierstrass, Abel--Poisson and Poisson wavelets and Mexican needlets. Boundedness of the uncertainty constant is in general not given, but it is a property of some of the wavelets from this class.
\end{abstract}

\begin{bfseries}Key Words and Phrases:\end{bfseries}   uncertainty product, time--frequency localization, spherical wavelets\\
\begin{bfseries}2010 MSC:\end{bfseries} 42C40

\section{Introduction}

Similarly as Heisenberg's uncertainty principle in quantum physics, several uncertainty principles are valid in mathematics.  They state that a function cannot be sharp both in space and in frequency, and quantitatively they are expressed by boundedness from below of an uncertainty constant whose definition depends on the domain of the function.

The problem of minimization of the uncertainty constant of wavelets on the real line -- whose analog on the sphere is the uncertainty product -- was studied since 1990s. The reason was that the wavelet transform gives information about time-frequency localization of continuous time-signals with finite energy and therefore the size of the time-frequency windows must be small. One of the first questions investigated in this context was to find wavelets preserving time-frequency localization as smoothness grows, where smoothness was understood for instance as H\"older exponent, as the decay rate of Fourier transform or as the number of vanishing moments of the wavelet or the scaling function, cf. \cite{eL11} and the references therein. A similar problem for periodic wavelets was considered in~\cite{LP14}. In my previous paper~\cite{IIN16UPW} I proved that the uncertainty constants of families of rotation-invariant Poisson wavelets over the $n$-dimensional sphere tend to the minimum value when some parameters approach their limits.

The present research is devoted to a broader class of spherical wavelets, and it is shown that even boundedness of the uncertainty constant is an exceptional property. I investigate functions given by
\begin{equation}\label{eq:wavelets}
\Psi_\rho(x)=\sum_{l=0}^\infty\left[\rho^a q_\nu(l)\right]^c\,e^{-\rho^a q_\nu(l)}\cdot\frac{l+\lambda}{\lambda}\,\mathcal C_l^\lambda(\cos\vartheta),
\end{equation}
where $x=(x_1,x_2,\dots,x_{n+1})\in\mathcal S^n\subseteq\mathbb R^{n+1}$ and~$\vartheta$ is the first hyperspherical variable, i.e.,
\begin{align*}
x_1&=\cos\vartheta,\\
x_2&=\sin\vartheta\cos\vartheta_2,\\
x_3&=\sin\vartheta\sin\vartheta_2\cos\vartheta_3,\\
&\dots\\
x_{n-1}&=\sin\vartheta\sin\vartheta_2\dots\sin\vartheta_{n-2}\cos\vartheta_{n-1},\\
x_n&=\sin\vartheta\sin\vartheta_2\dots\sin\vartheta_{n-2}\sin\vartheta_{n-1}\cos\varphi,\\
x_{n+1}&=\sin\vartheta\sin\vartheta_2\dots\sin\vartheta_{n-2}\sin\vartheta_{n-1}\sin\varphi.
\end{align*}
Functions depending only on~$\vartheta$ are invariant with respect to a rotation around the $X_1$-axis and they are called \emph{zonal} (or \emph{rotation-invariant}). In formula~\eqref{eq:wavelets}, $\rho\in\mathbb R_+$ , $q_\nu$ is a polynomial of degree~$\nu$, strictly positive and monotonously increasing for $l\geq1$, and $a$, $c$ --- some positive constants. Further, $\lambda$ a parameter depending on space dimension, $\lambda=\frac{n-1}{2}$, and $\mathcal C_l^\lambda$ denotes Gegenbauer polynomial of degree~$l$ and order~$\lambda$.

There are two essentially different definitions of continuous spherical wavelets, namely those based on group theoretical approach~\cite{AVn} and those derived from approximate identities~\cite{IIN15WT}, see discussion in~\cite[Sec.~5]{IIN15WT}. It is shown in \cite[Theorem~2.6]{IIN16FDW} that a wide class of (non-zonal) wavelets belonging to the latter class can be derived from functions defined by~\eqref{eq:wavelets}. It can be easily verified that the zonal functions~\eqref{eq:wavelets} are wavelets as well. Although \cite[Theorem~2.6]{IIN16FDW} concerns a wider class of functions, we concentrate on those with a polynomial in the exponent of the exponential function. The reason is, the wavelets used in applications and considered in the literature, such as Gauss--Weierstrass wavelets, Abel--Poisson wavelets or the whole family of Poisson wavelets \cite{FGS-book,HI07,IIN15PW}, as well as Mexican needlets \cite{GM09a,GM09b} are of this type. Further, the aim of the present research is to compare wavelet families, characterized by parameters $a$, $c$, $\nu$, and the polynomial~$q$, in respect of their uncertainty product. Zonal wavelets are in the author's opinion a class of functions wide enough to yield interesting results. Computation of the uncertainty product of directional wavelets will be the object of my future research.

An uncertainty principle for twice continuously differentiable functions on the two-dimensional sphere was derived by Narcowich and Ward in~\cite{NW96}. The definitions of variances in space and momentum domain (whose square roots are \emph{uncertainties}) were used by other authors to characterize distinct function families in respect of their localization, see e.g. \cite{FM99}. The case of $n$-dimensional spheres was investigated in~\cite{RV97}, where an uncertainty principle for zonal functions was derived. These ideas were further generalized in \cite{GG04a,GG04b} to continuously differentiable non-vanishing functions over~$\mathcal S^n$. The uncertainty principle stated in \cite{NW96,RV97} is weaker than that from \cite{GG04a,GG04b}. To my best knowledge, the notion \emph{uncertainty product} was introduced in \cite{nLF03}. Contrary to the authors of \cite{NW96,RV97}, La\'in Fern\'andez uses the word \emph{variance} for what was called \emph{uncertainty} and the product of variances in space and in momentum domain is called the \emph{uncertainty product} of a twice continuously differentiable function. Except for these distinct expressions, the definitions coincide in the above mentioned papers.

A series representation of the uncertainty product using Gegenbauer coefficients of a zonal function over $\mathcal S^2$ is derived in \cite{nLF03}. Similarly, a series representation of the uncertainty product  of a function over $\mathcal S^1$ is presented in~\cite{GG04a,GG04b} and for zonal  $\mathcal S^n$-functions in \cite{IIN16MR}.

\section{Preliminaries}\label{sec:sphere}

Let $\mathcal{S}^n$ denote the $n$--dimensional unit sphere in $(n+1)$--dimensional Euclidean space~$\mathbb{R}^{n+1}$ and $\Sigma_n=\frac{2\pi^{\lambda+1}}{\Gamma(\lambda+1)}$ is its Lebesgue measure. Gegenbauer polynomials $C_l^\lambda$ of order~$\lambda\in\mathbb R$, and degree~$l\in\mathbb{N}_0$, are defined in terms of their generating function
$$
\sum_{l=0}^\infty C_l^\lambda(t)\,r^l=\frac{1}{(1-2tr+r^2)^\lambda},\qquad t\in[-1,1].
$$
They are real-valued and for a fixed $\lambda\ne0$ orthogonal to each other with respect to the weight function~$t\mapsto\left(1-t^2\right)^{\lambda-\frac{1}{2}}$.
Integrable zonal functions over the sphere have the Gegenbauer expansion
$$
f(\cos\vartheta)=\sum_{l=0}^\infty\widehat f(l)\,C_l^\lambda(\cos\vartheta)
$$
with Gegenbauer coefficients
\begin{equation}\label{eq:Gegenbauer_coeffs}
\widehat f(l)=c(l,\lambda)\int_{-1}^1 f(t)\,C_l^\lambda(t)\left(1-t^2\right)^{\lambda-1/2}dt,
\end{equation}
where~$c$ is a constant that depends on~$l$ and~$\lambda=\frac{n-1}{2}$.

The variances in space and momentum domain of a $\mathcal C^2(\mathcal S^n)$--function~$f$ with $\int_{\mathcal S^n}x\,|f(x)|^2\,d\sigma(x)\ne0$ are given by \cite{nLF03}
$$
\text{var}_S(f)=\left(\frac{\int_{\mathcal S^n}|f(x)|^2\,d\sigma(x)}{\int_{\mathcal S^n}x\,|f(x)|^2\,d\sigma(x)}\right)^2-1
$$
and
$$
\text{var}_M(f)=-\frac{\int_{\mathcal S^n}\Delta^\ast f(x)\cdot \bar f(x)\,d\sigma(x)}{\int_{\mathcal S^n}|f(x)|^2\,d\sigma(x)},
$$
where $\Delta^\ast$ is the Laplace--Beltrami operator on~$\mathcal S^n$. The quantity
$$
U(f)=\sqrt{\text{var}_S(f)}\cdot\sqrt{\text{var}_M(f)}
$$
is called the uncertainty product of~$f$.

The uncertainty product of zonal functions may be computed from their Gegenbauer coefficients \cite[Lemma~4.2]{IIN16MR} and according to the spherical uncertainty principle it is bounded from below by $\frac{n}{2}$, see \cite{NW96,RV97} and \cite[formula (4.37)]{GG04a}, \cite[formula (12)]{GG04b}.

\begin{lem}\label{lem:varS_varM} Let a zonal  square integrable and continuously differentiable function over $\mathcal S^n$ be given by its Gegenbauer expansion
$$
f(\cos\vartheta)=\sum_{l=0}^\infty\widehat f(l)\,\mathcal C_l^\lambda(\cos\vartheta).
$$
Its variances in space and momentum domain are equal to
\begin{align}
\text{var}_S(f)&=\left(\frac{\sum_{l=0}^\infty\frac{\lambda}{l+\lambda}\,\binom{l+2\lambda-1}{l}\,|\widehat f(l)|^2}{\sum_{l=0}^\infty\binom{l+2\lambda}{l}\,
   \frac{\lambda^2\left[\overline{\widehat f(l)}\,\widehat f(l+1)+\widehat f(l)\,\overline{\widehat f(l+1)}\right]}{(l+\lambda)(l+\lambda+1)}}\right)^2-1,\label{eq:varS}\\
\text{var}_M(f)&=\frac{\sum_{l=1}^\infty\frac{ l\lambda(l+2\lambda)}{l+\lambda}\,\binom{l+2\lambda-1}{l}\,|\widehat f(l)|^2}
   {\sum_{l=0}^\infty\frac{\lambda}{l+\lambda}\,\binom{l+2\lambda-1}{l}\,|\widehat f(l)|^2},\label{eq:varM}
\end{align}
whenever the series are convergent.
\end{lem}

\begin{thm}For $f\in\mathcal L^2(\mathcal S^n)\cap\mathcal C^1(\mathcal S^n)$, $U(f)\geq\frac{n}{2}$.
\end{thm}

\section{The uncertainty product of spherical wavelets}\label{sec:uncertainty}

For the proof of the main statement of this paper, Theorem~\ref{thm:UCwv}, we need some lemmas.

\begin{lem}\label{lem:asymptotics_f_powers}Let $\nu$, $d$, and $\rho$ be positive numbers. Then for each $k\in\mathbb N$
\begin{equation}\label{eq:asymptotics}
\sum_{l=k}^\infty l^de^{-\rho l^\nu}=\frac{\Gamma\left(\tfrac{d+1}{\nu}\right)}{\nu}\,\rho^{-\tfrac{d+1}{\nu}}+\mathcal O\left(\rho^{-\tfrac{d}{\nu}}\right)\qquad\text{for }\rho\to0.
\end{equation}
\end{lem}
\begin{bfseries}Proof.\end{bfseries} Consider the function
\begin{align*}
f:\,\mathbb R_+&\to\mathbb R_+,\\
x&\mapsto x^de^{-\rho x^\nu}.
\end{align*}
It is monotonously increasing for $x\in\left(0,\sqrt[\nu]{\frac{d}{\rho \nu}}\right)$ and monotonously decreasing for $x\in\left(\sqrt[\nu]{\frac{d}{\rho \nu}},+\infty\right)$. Consequently,
\begin{equation}\label{eq:estimation_quadrature}
\int_0^\infty f(x)\,dx-(k+1)\cdot f_\text{max}\leq\sum_{l=k}^\infty f(l)\leq\int_0^\infty f(x)\,dx+f_\text{max}.
\end{equation}
Since
\begin{equation}\label{eq:integral_xdexp}
\int_0^\infty f(x)\,dx=\frac{\Gamma\left(\tfrac{d+1}{\nu}\right)}{\nu}\,\rho^{-\tfrac{d+1}{\nu}}
\end{equation}
and
\begin{equation}\label{eq:fmax}
f_\text{max}=f\left(\sqrt[\nu]{\frac{d}{\rho \nu}}\right)=\left(\frac{d}{\nu e}\right)^\frac{d}{\nu}\rho^{-\frac{d}{\nu}},
\end{equation}
one obtains~\eqref{eq:asymptotics}.\hfill$\Box$

\begin{lem}\label{lem:asymptotics_f_power+polynomial}Let $d$ and $\rho$ be positive numbers, and $q(l)=l^\nu+a_{\nu-1}l^{\nu-1}+\dots+a_0$ --  a polynomial of degree $\nu\geq1$, positive and monotonously increasing for $l\geq1$. Then for each $k\in\mathbb N$
\begin{equation}\label{eq:asymptotics_q}
\sum_{l=k}^\infty l^de^{-\rho q(l)}=\frac{\Gamma\left(\tfrac{d+1}{\nu}\right)}{\nu}\,\rho^{-\tfrac{d+1}{\nu}}+\mathcal O\left(\rho^{-\tfrac{d}{\nu}}\right)\qquad\text{for }\rho\to0.
\end{equation}
\end{lem}

\begin{bfseries}Proof. \end{bfseries}If $\nu\geq2$, consider the function
\begin{align*}
f:\,\mathbb R_+&\to\mathbb R_+,\\
x&\mapsto x^de^{-\rho q(x)}.
\end{align*}
Its derivative
$$
f^\prime(x)=x^{d-1}\,e^{-\rho q(x)}\left[d-x\rho q^\prime(x)\right]
$$
is positive for $x=0$, tends to $-\infty$ for $x\to+\infty$ and changes the sign in~$x_0$ such that $x_0 q^\prime(x_0)=\frac{d}{\rho}$. Thus, the function~$f$ has a local maximum in $x_0$. Consequently, estimation~\eqref{eq:estimation_quadrature} holds also in this case.

For $\rho\in(0,1)$ set $$R=R(\rho)=\rho^{\frac{1}{2(1-\nu)}}.$$ Then
$$
|q(x)-x^\nu|\leq aR^{\nu-1}=\frac{a}{\sqrt\rho}
$$
for $x\in(0,R)$ and $a=\sum_{\iota=0}^{\nu-1}|a_\iota|$. Further, for each $\epsilon>0$ there exists $\rho_0\in(0,1)$ such that
$$
\left|\frac{q(x)}{x^\nu}-1\right|\leq\epsilon
$$
 for all $\rho\in(0,\rho_0)$ and $x\in[R(\rho),\infty)$. Consequently,
\begin{align*}
&x^de^{-\rho x^\nu}\,e^{-a\sqrt\rho}\leq f(x)\leq x^de^{-\rho x^\nu}\,e^{a\sqrt\rho}\qquad\text{for }x\in(0,R)\text{ and}\\
&x^d\left(e^{-\rho x^\nu}\right)^{1+\epsilon}\leq f(x)\leq x^d\left(e^{-\rho x^\nu}\right)^{1-\epsilon}\qquad\text{for }x\in[R,\infty),
\end{align*}
and further
\begin{equation}\label{eq:estimations_integrals}\begin{split}
e^{-a\sqrt\rho}\int_0^R x^de^{-\rho x^\nu}\,dx&\leq\int_0^R f(x)\,dx\leq e^{a\sqrt\rho}\int_0^R x^de^{-\rho x^\nu}\,dx,\qquad\text{and}\\
\int_R^\infty x^de^{-(1+\epsilon)\rho x^\nu}\,dx&\leq\int_R^\infty f(x)\,dx\leq\int_R^\infty x^de^{-(1-\epsilon)\rho x^\nu}\,dx.
\end{split}\end{equation}
The sum of the left-hand-sides of~\eqref{eq:estimations_integrals} can be estimated from below as follows:
\begin{align*}
e^{-a\sqrt\rho}&\int_0^R x^de^{-\rho x^\nu}\,dx+\int_R^\infty x^de^{-(1+\epsilon)\rho x^\nu}\,dx\\
&\geq e^{-a\sqrt\rho}\int_0^R x^de^{-(1+\epsilon)\rho x^\nu}\,dx+e^{-a\sqrt\rho}\int_R^\infty x^de^{-(1+\epsilon)\rho x^\nu}\,dx\\
&=e^{-a\sqrt\rho}\int_0^\infty x^de^{-(1+\epsilon)\rho x^\nu}\,dx,
\end{align*}
and for the right-hand-side of~\eqref{eq:estimations_integrals} we obtain in an analogous way
$$
e^{a\sqrt\rho}\int_0^R x^de^{-\rho x^\nu}\,dx+\int_R^\infty x^de^{-(1-\epsilon)\rho x^\nu}\,dx\leq e^{a\sqrt\rho}\int_0^\infty x^de^{-(1-\epsilon)\rho x^\nu}\,dx.
$$
Consequently, by~\eqref{eq:integral_xdexp},
$$
e^{-a\sqrt\rho}\,\frac{\Gamma\left(\frac{d+1}{\nu}\right)}{\nu}\,[(1+\epsilon)\rho]^{-\frac{d+1}{\nu}}\leq\int_0^\infty f(x)\,dx
   \leq e^{a\sqrt\rho}\,\frac{\Gamma\left(\frac{d+1}{\nu}\right)}{\nu}\,[(1-\epsilon)\rho]^{-\frac{d+1}{\nu}}.
$$
Since we investigate the behavior of~$\sum_k^\infty f(l)$ for $\rho\to0$, we can assume without loss of generality that $\epsilon\to0$. Thus
$$
\int_0^\infty x^de^{-\rho q(x)}\,dx\to\frac{\Gamma\left(\frac{d+1}{\nu}\right)}{\nu}\,\rho^{-\frac{d+1}{\nu}}\qquad\text{for }\rho\to0.
$$
In order to estimate~$f_\text{max}$ note that since $q(x)>0$ for $x\geq1$, there exists $\alpha>0$ such that
$$
q(x)\geq\alpha x^\nu\qquad\text{for }x\geq1.
$$
Thus,
$$
f(x)\leq x^de^{-\alpha\rho x^\nu}.
$$
Apply~\eqref{eq:fmax} to the function on the right-hand-side of this inequality to get
$$
f(x)\leq\left(\frac{d}{\alpha \nu e}\right)^\frac{d}{\nu}\,\rho^{-\frac{d}{\nu}}.
$$
We obtain~\eqref{eq:asymptotics_q} by the same arguments as in the proof of Lemma~\ref{lem:asymptotics_f_powers}.

For $\nu=1$ write the series as
$$
\sum_{l=k}^\infty l^de^{-\rho(l+a_0)}=e^{-\rho a_0}\sum_{l=k}^\infty l^de^{-\rho l}
$$
and apply Lemma~\ref{lem:asymptotics_f_powers}.
\hfill$\Box$

\begin{lem}\label{lem:asymptotics_f_irrational+polynomial}Let $d$ and $\rho$ be positive numbers, $q(l)=l^\nu+a_{\nu-1}l^{\nu-1}+\dots+a_0$ --  a polynomial of degree $\nu\geq1$, positive and monotonously increasing for $l\geq1$, $Q(l)=l^r+b_{r-1}l^{r-1}+\dots+b_0$ -- a polynomial of degree~$r\geq1$, positive for positive~$l$, and $p\in\mathbb N$. Then
\begin{equation*}
\sum_{l=1}^\infty l^pQ(l)^de^{-\rho q(l)}
   =\frac{\Gamma\left(\tfrac{p+rd+1}{\nu}\right)}{\nu}\,\rho^{-\tfrac{p+rd+1}{\nu}}+\mathcal O\left(\rho^{-\frac{p+rd}{\nu}}\right)\qquad\text{for }\rho\to0.
\end{equation*}
\end{lem}

\begin{bfseries}Proof.\end{bfseries} Set
$$
k=\left[\max_{\iota\in\{r-1,r-2,\dots,1,0\}}|b_\iota|\right]+1,
$$
$[\circ]$ denoting the \emph{entier} function.
Then, for $l\geq1$,
$$
Q(l)\leq l^r+rkl^{r-1}\leq(l+k)^r
$$
and
$$
Q(l)\geq l^r-rkl^{r-1}.
$$
Note that for positive $l$ and $c$, and $r\geq2$
\begin{align*}
(l-c)^r&\leq l^r-crl^{r-1}+2^r\sum_{\iota=2}^rc^\iota l^{r-\iota}\\
   &\leq l^r-crl^{r-1}+2^r\cdot\frac{c^2l^{r-2}}{1-\frac{c}{l}}\\
   &\leq l^r-\frac{cr}{2}\,l^{r-1}\qquad\text{for }l\geq\left(\frac{2^{r+1}}{r}+1\right)c,
\end{align*}
and for $r=1$
$$
l-c\leq l-\frac{c}{2}.
$$
Consequently, for $l\geq2k\left(\frac{2^{r+1}}{r}+1\right)$,
$$
Q(l)\geq\left(l-2k\right)^r.
$$
Thus, by Lemma~\ref{lem:asymptotics_f_power+polynomial},
\begin{align}
\sum_{l=1}^\infty &l^pQ(l)^de^{-\rho q(l)}\leq\sum_{l=1}^\infty(l+k)^{p+rd}\,e^{-\rho q(l)}=\sum_{l=1+k}^\infty l^{p+rd}\,e^{-\rho q(l-k)}\notag\\
&\leq\frac{\Gamma\left(\frac{p+rd+1}{\nu}\right)}{\nu}\,\rho^{-\frac{p+rd+1}{\nu}}+\mathcal O\left(\rho^{-\frac{p+rd}{\nu}}\right)\qquad\text{for }\rho\to0.\label{eq:ub_est_poly_poly}
\end{align}
On the other hand, set $K=\left[2k\left(\frac{2^{r+1}}{r}+1\right)\right]+1$. Then,
\begin{align}
\sum_{l=1}^\infty &l^pQ(l)^de^{-\rho q(l)}\geq\sum_{l=1+K}^\infty l^pQ(l)^de^{-\rho q(l)}\geq\sum_{l=1+K}^\infty(l-2k)^{p+rd}\,e^{-\rho q(l)}\notag\\
&=\sum_{l=1+K-2k}^\infty l^{p+rd}\,e^{-\rho q(l+2k)}
   \geq\frac{\Gamma\left(\frac{p+rd+1}{\nu}\right)}{\nu}\,\rho^{-\frac{p+rd+1}{\nu}}
   +\mathcal O\left(\rho^{-\frac{p+rd}{\nu}}\right)\qquad\text{for }\rho\to0.\label{eq:lb_est_poly_poly}
\end{align}
One obtains the assertion from \eqref{eq:ub_est_poly_poly} and \eqref{eq:lb_est_poly_poly}.\hfill$\Box$

After this preparation the main result may be proved.

\begin{thm}\label{thm:UCwv} Let $\{\Psi_\rho\}$ be a zonal wavelet family with
\begin{equation}\label{eq:Gegenbauer_coeffs}
\widehat\Psi_\rho(l)=\left[\rho^a q_\nu(l)\right]^c\,e^{-\rho^a q_\nu(l)}\cdot\frac{l+\lambda}{\lambda},
\end{equation}
where $a>0$, $c>0$, and $q_{\nu}(l)=a_\nu l^\nu+a_{\nu-1}l^{\nu-1}+\dots+a_1 l+a_0$ is a polynomial of degree~$\nu$, positive and monotonously increasing for $l\geq1$. The
uncertainty product of~$\Psi_\rho$ for $\rho\to0$ behaves like
$$
U(\Psi_\rho)\leq\mathcal O\left(\rho^{\frac{-a}{2\nu}}\right).
$$
\end{thm}

\begin{bfseries}Proof. \end{bfseries}Note that in the expressions \eqref{eq:varS} and \eqref{eq:varM}, in the numerators and the denominators there occur products of two Gegenbauer coefficients of~$f$ such that the factor~$\frac{\rho^{ac}}{\lambda}$ in~\eqref{eq:Gegenbauer_coeffs} can be neglected by computation of~var$_S(\Psi_\rho)$ and~var$_M(\Psi_\rho)$. Consider
$$
\widehat f_\rho(l)=(l+\lambda)\,q_\nu(l)^c\,e^{-\rho^a q_\nu(l)}
$$
and compute the numerator of~var$_S(f_\rho)$,
$$
N_S(f_\rho):=\sum_{l=0}^\infty\frac{\lambda}{l+\lambda}\binom{l+2\lambda-1}{l}|\widehat f_\rho(l)|^2
   =\lambda\sum_{l=0}^\infty P_{2\lambda-1}(l)\,(l+\lambda)\,q_{\nu}(l)^{2c}\,e^{-2\rho^aq_\nu(l)},
$$
where~$P_{2\lambda-1}$ is a polynomial of degree~$2\lambda-1$ and with leading coefficient equal to~$\frac{1}{(2\lambda-1)!}$. According to Lemma~\ref{lem:asymptotics_f_irrational+polynomial}, for $\rho\to0$ it behaves like
$$
N_S(f_\rho)=\frac{\lambda\,a_\nu^{2c}\,\Gamma\left(\frac{2(\lambda+\nu c)+1}{\nu}\right)}{(2\lambda-1)!\,\nu}\,(2\rho^aa_\nu)^{-\frac{2(\lambda+\nu c)+1}{\nu}}
   +\mathcal O\left(\rho^{-\frac{2a(\lambda+\nu c)}{\nu}}\right).
$$
In a similar manner, for the denominator of~\eqref{eq:varS} we obtain
\begin{align*}
D_S(f_\rho)&:=\sum_{l=0}^\infty\binom{l+2\lambda}{l}\frac{2\lambda^2\widehat f(l)\widehat f(l+1)}{(l+\lambda)(l+\lambda+1)}\\
&=2\lambda^2\sum_{l=0}^\infty P_{2\lambda}(l)\left[q_\nu(l)\,q_\nu(l+1)\right]^c\,e^{-\rho^a[q_\nu(l)+q_\nu(l+1)]}
\end{align*}
with a polynomial $P_{2\lambda}$ of degree~$2\lambda$ and with leading coefficient~$\frac{1}{(2\lambda)!}$. Consequently,
$$
D_S(f_\rho)=\frac{\lambda\,a_\nu^{2c}\,\Gamma\left(\frac{2(\lambda+\nu c)+1}{\nu}\right)}{(2\lambda-1)!\,\nu}\,(2\rho^aa_\nu)^{-\frac{2(\lambda+\nu c)+1}{\nu}}
   +\mathcal O\left(\rho^{-\frac{2a(\lambda+\nu c)}{\nu}}\right).
$$
Thus,
\begin{equation}\label{eq:asymptotics_varS}
\text{var}_S(f_\rho)\leq\left[1+\mathcal O\left(\rho^{\frac{a}{\nu}}\right)\right]^2-1=\mathcal O\left(\rho^{\frac{a}{\nu}}\right)
\qquad\text{for }\rho\to0.
\end{equation}
Further, for the numerator of~$\text{var}_M(f_\rho)$ we have
\begin{align*}
N_M(f_\rho)&:=\sum_{l=1}^\infty\frac{l\lambda(l+2\lambda)}{l+\lambda}\binom{l+2\lambda-1}{l}|\widehat f(l)|^2\\
&=\lambda\sum_{l=1}^\infty l(l+\lambda)(l+2\lambda)P_{2\lambda-1}(l)\,q_\nu(l)^{2c}\,e^{-2\rho^aq_\nu(l)}\\
&=\frac{\lambda\,a_\nu^{2c}\,\Gamma\left(\frac{2(\lambda+\nu c)+3}{\nu}\right)}{(2\lambda-1)!\,\nu}\,(2\rho^aa_\nu)^{-\frac{2(\lambda+\nu c)+3}{\nu}}
   +\mathcal O\left(\rho^{-\frac{2a(\lambda+1+\nu c)}{\nu}}\right),
\end{align*}
and the denominator of~var$_M(f_\rho)$ is equal to~$N_S$,
$$
D_M(f_\rho)=\frac{\lambda\,a_\nu^{2c}\,\Gamma\left(\frac{2(\lambda+\nu c)+1}{\nu}\right)}{(2\lambda-1)!\,\nu}\,(2\rho^aa_\nu)^{-\frac{2(\lambda+\nu c)+1}{\nu}}
   +\mathcal O\left(\rho^{-\frac{2a(\lambda+\nu c)}{\nu}}\right).
$$
Consequently,
\begin{equation}\label{eq:asymptotics_varM}
\text{var}_M(f_\rho)\leq\frac{\Gamma\left(\frac{2(\lambda+\nu c)+3}{\nu}\right)}{\Gamma\left(\frac{2(\lambda+\nu c)+1}{\nu}\right)}\,
   (2\rho^aa_\nu)^{-\frac{2}{\nu}}+\mathcal O\left(1\right).
\end{equation}
It follows from~\eqref{eq:asymptotics_varS}, \eqref{eq:asymptotics_varM}, and the definition of the uncertainty product that
$$
U(\Psi_\rho)=U(f_\rho)\leq\mathcal O\left(\rho^{\frac{-a}{2\nu}}\right)\qquad\text{for }\rho\to0.
$$
\hfill$\Box$

\section{Discussion}
It is apparent that the uncertainty product of a wavelet family given by \eqref{eq:Gegenbauer_coeffs} is in general unbounded for $\rho\to0$. Computation of the  second terms in the expansions of~$N_S(f_\rho)$ and $D_S(f_\rho)$ (in order to  obtain the exact coefficient in the first term of the representation \eqref{eq:varS}) is impossible without explicit knowledge of the coefficients of the polynomial~$q$, and even if they were known, it is quite sophisticated, see \cite{LFP03}. Note that in the case of Poisson wavelets
$$
g_\rho^m(x)=\frac{1}{\Sigma_n}\sum_{l=0}^\infty\frac{l+\lambda}{\lambda}\,(\rho l)^me^{-\rho l}\mathcal C_l^\lambda(\cos\vartheta),
$$
the term $\rho^{\frac{a}{\nu}}=\rho$ disappears, and one has var$_S(g_\rho^m)=\mathcal O(\rho^2)$, see~\cite{IIN16UPW}. Consequently, Poisson wavelets have bounded uncertainty product for $\rho\to0$.

Similarly, Gauss--Weierstrass kernel on the two-dimensional sphere
$$
\Phi_\rho(x)=\frac{1}{4\pi}\sum_{l=0}^\infty(2l+1)\, e^{-\rho l(l+1)}\,\mathcal C_l^{1/2}(\cos\vartheta)
$$
has a bounded uncertainty product, see \cite{LFP03}.

Further, the exponent~$c$ may have an influence on the value of the coefficients in the exact expansion of the variances and the uncertainty product, as it is the case by Poisson wavelets, but the exponents in these expansions are dependent only on~$a$ and~$\nu$. However, the constant~$a$ can be easily ousted: Set $\widetilde\rho:=\rho^a$. Since the measures $\alpha(\widetilde\rho)$ and $\alpha(\rho)$ differ only by a multiplicative constant, $\rho$ may be replaced by~$\widetilde\rho$ in the investigation of the uncertainty product.

In the general case, the increase of the uncertainty product of a wavelet for $\rho\to0$ is undesirable. Thus, bounded uncertainty product of Poisson wavelets \cite{HI07,IIN15PW} is the next property (beside explicit expressions in terms of spherical variables $\vartheta$ and $\varphi$ \cite{IIN15PW}, Euclidean limit property \cite{IIN15WT}, and existence of discrete frames with density proportional to the scale \cite{IH10,IIN16WF}) that makes this wavelet families superior to other ones. On the other hand, the uncertainty product of another wavelet families must be computed. Note that the computations in~\cite{IIN16UPW} cannot be applied to Abel--Poisson wavelet (that can be regarded as Poisson wavelet of order \mbox{$m=\frac{1}{2}$}). Similarly, to the author's best knowledge, the uncertainty product of Gauss--Weierstrass wavelet
$$
\Psi_\rho(x)=\frac{1}{4\pi}\sum_{l=0}^\infty(2l+1)\,\sqrt{2\rho l(l+1)} e^{-\rho l(l+1)}\,\mathcal C_l^{1/2}(\cos\vartheta)
$$
has not been computed so far (contrary to the uncertainty product of Gauss--Weierstrass \emph{kernel}). Since these wavelet families are most popular, it seems to be important to characterize them in respect of their uncertainty product.

\end{document}